\documentstyle{amsppt} 
\magnification\magstep1

\TagsOnRight
\NoRunningHeads

\catcode`\@=11
\def\logo@{}
\let\varindent@\indent
\catcode`\@=13

\def\ff{\mathrel{=\!\!\!\mid\,}} 
\def\ffm#1{\mathrel{=_{\scriptscriptstyle #1}\hskip-4.5mm\mid\hskip 3mm}}

\def\Bo{{\mathaccent'27B}}
\def\Ko{{\mathaccent'27K}}
\def\bR{{\Bbb R}}

\def\card{\operatorname{card}}
\def\const{\operatorname{const}}
\def\cross{\operatorname{cr}}
\def\eqc{\operatorname{eqc}}
 \let\le\leqslant
 \let\ge\geqslant
 \let\eps\varepsilon
\let\phi\varphi
\let\<\langle  \let\>\rangle

\topmatter
\title{On strict suns in $\ell^\infty(3)$}
\endtitle
\author{A.\,R.\,Alimov}
\thanks Supported by the Russian Fund for Basic Researches,
Project No.~02-01-00248.
\endthanks
\endauthor
\keywords Strict sun; sun; $\ell^\infty(n)$; characterisation; best approximation;
metrical convexity; $\ell^1$-convexity
\endkeywords
\subjclass 41A65 \endsubjclass
\abstract\nofrills\indent
A subset $M$ of a normed linear space $X$ is said to be a 
{\it strict sun\/} if, for every point $x\in X\setminus M$, 
the set of its nearest points from~$M$ is non-empty and 
if $y\in M$ is a nearest point from~$M$ to ~$x$, then $y$
is a nearest point from~$M$ to all points from the ray $\{\lambda x+(1-
\lambda)y\,|\, \lambda>0\}$. In the paper there obtained a geometrical
characterisation of strict suns in  $\ell^\infty(3)$. In comparison 
with [1] we establish a more precise property of stict suns. 
\endabstract
\endtopmatter
\document

A subset $M$ of a normed linear space~$X$ is said to be a 
{\it strict sun} if, for every $x\in 
X\setminus M$, the set of its nearest points from~$M$ is 
non-empty and if $y\in M$ is the nearest point to~$x$, then $y$ is 
the nearest ponit from~$M$ to every point from the ray starting 
at~$y$ and passing through~$x$.
In this paper we obtain a geometric characterisation of 
strict suns in the space $\ell^\infty(3)$ 
(Theorem~1). This characterisation provides a more precise
property of stricts suns in $\ell^\infty(3)$ than Theorem~A (see 
below) from~[1] does.
\medskip

\subhead{1. Definitions and notation}\endsubhead
We will consider only real spaces.

Following K.~Menger [2] and  H.~Berens and L.~Hetzelt [3],
a set $M\subset\bR^n$ is called $\ell^1$-{\it 
convex} if, for all $x,y\in M$, $x\ne y$, there is a point 
$z\in M$, $z\ne x$, $z\ne y$, such that  $\|x-y\|_{\ell^1}=\|x-
z\|_{\ell^1}+\|z-y\|_{\ell^1}$ (here $\|{\cdot}\|_{\ell^1}$ is the  
standard $\ell^1$-norm on~$\bR^n$).

As it is shown in [3], $\ell^1$-convexity (together with closedness) proved 
to be a characteristic property of suns~$\ell^\infty(n)$. (Here  
we recall that a  set $M\subset X$ is a {\it sun} if, for every point  
$x\in X\setminus M$ there is a point $y$, that is nearest from~$M$ 
to~$x$ and such that $y$ is a nearest point to all points 
from the  ray starting at~$y$ and passing through~$x$;
such a point $y$ is called a {\it solar} point for~$x$). It is clear that
every strict sun is a sun.  The inverse statement is not true 
in general. On suns, strict suns and other approximative sets 
see [4],~[5]. 

An elegant characterisation of suns in $\ell^\infty(n)$ was obtained by
H.~Berens and L.~Hetzelt~[3] (see also [6]). 

\proclaim{Theorem \smc (H.~Berens, L.~Hetzelt)} A closed nonvoid set 
$M\subset\bR^n$ is a sun in~$\ell^\infty(n)$ if and only if it is 
$\ell^1$-convex. 
\endproclaim

The first characterisation of strict suns in $\ell^\infty(n)$ in geometrical
terms was independently obtained by B.~Brosowski~[7] and Ch.~Dunham~[8]. 
In~[1] (Theorem~A) there obtained a characterisation of strict suns  
in~$\ell^\infty(n)$, which is similar to the metric characterisation due
to Berens and Hetzelt.

Here some necessary definitions are required. 

Given a set $M\subset\bR^n$, we define 
$$
  \gather
        \eqc(M)=\eqc_{\bR^n}(M)
        :=\big\{j\in\{1,\dots,n\}\,|\, x_j=y_j\ 
        \text{for\ all}\ x,y\in M \big\}, \\
        x=(x_1,\dots,x_n),\ \  y=(y_1,\dots,y_n);
  \endgather
$$
the abbreviation ``eqc'' stands for  `equal coordinates'.

For $x,y\in\bR^n$ and $M\subset\bR^n$, we define
\roster
  \item""$x\ffm M y$, if $x_j\ne y_j$ for all $j\notin\eqc(M)$;
  \item""$x\ff y$,    if $x_i\ne y_i$ for all $i=1,\dots,n$.
\endroster

A set $M$ is said to be  {\it strictly $\ell^1$-convex} if
\roster
 \item"{1)}"~it is $\ell^1$-convex, and 
 \item"{2)}"{for every $x,y\in M$, $x\ffm M y$, there is an 
$\ell^1$-geodesic segment $k(t)\subset M$, $0\le t\le 1$, connecting
$x$ and~$y$ for which the functions $[0,1]\ni t\mapsto k_j(t)$ are monotonic
on~$t$ for $j=1,2,\dots,n$ and strictly monotonic for $j\notin \eqc(M)$.} 
\endroster

\noindent

In our terminology, the $\ell^1$-convexity of~$M$ means~[3] 
that any two points  $x,y\in M$ can be connected by an $\ell^1$-geodesic
segment $k(\cdot)\subset M$, with monotonic in~$t$ coordinate functions 
$k_i(t)$ (which are not necessarily strictly monotonic).

Following [1], by a {\it main cocross}\/ in~$\bR^n$ we will understand
the set
$$
        c^n(x)=\big\{y\in\bR^n\,|\,
        \card\eqc(y,x)\ge  1\}, \ \ x\in \bR^n \ \text{is\ fixed},
$$
and by a  {\it cocross $c_{J}(x)$ with respect to the assembly of coordinates} 
$J=\{j_1,\dots, j_k\} \subset\{1,\dots,n\}$  we will understand the set
$$
        c_J(x)=c^n_J(x)=\bigg\{y\in\bR^n\,\bigg|\,
        \aligned
                        & y_j=\const(j)\\
                        &\card(\eqc(x,y)\cap J)\ge 1\
        \endaligned\
        \aligned
                &\text{for}\ j\in J, \\
                &\text{for}\ i\notin J.
        \endaligned
                                \bigg\}
$$
In addition, to exclude degenerate cases, we will {\it 
always\/} assume that $k\le n-2$; i.e., $\card J\le n-2$ in the definition
of $c_{J}(x)$. Finally, a {\it cocross} is any set $C\subset\bR^n$ 
such that $C\subset c_{\eqc(C)}(x)$ for some  
$x\in\bR^n$. (Note that always $\card\eqc(C)<n-1$.)

For $x\in \bR^n$ we will define 
$$
    \cross(x)= \cross^n(x)=
    \{y\in\bR^n\,|\, \card \eqc(x,y)\ge n-1\}\ 
    \text{a\ {\it main\ cross}\ in}\ \bR^n.
$$
A {\it cross} $C$ is a subset of some main cross in~$\bR^n$
such that $\eqc(C)=\emptyset$; we will call  it {\it $\bR^n$-cross\/}.
It is clear that any $\bR^n$-cross is a cocross.  

\medskip

The next theorem [1] characterises strict suns in~$\ell^\infty(n)$.

\proclaim{Theorem A} Let $\emptyset\ne M\subset\bR^n$ be closed. 
Then $M$ is a strict sun in $\ell^\infty(n)$ if and only if it is 
{\rm 1)}~strictly $\ell^1$-convex and {\rm 2)}~not a cocross.
\endproclaim

Conditions 1) and 2) are not independent ones. In~[1] it was mentioned 
that $\bR^3$-cross $\cross^3(x)$  is an $\ell^1$-convex set. 
However, the main cocross $c^n(x)$ is, clearly, 
not $\ell^1$-convex. 

\subhead{2. The main theorem} \endsubhead

Here we prove a result in which the characteristical conditions in  
Theorem~A are rectified for $n\le 3$. Namely, condition~2) is replaced
with the more strong one: {\it a set is not an $\bR^3$-cross}. 
We also prove Proposition~1, which characterises strict 
$\ell^1$-convex cocrosses in~$\bR^3$ among all cocrosses in~$\bR^3$.

\proclaim{Theorem 1} Let $\emptyset\ne M\subset\bR^3$ be closed.
Then $M$ is a strict sun in~$\ell^\infty(3)$ if and only if 
it is {\rm 1)}~strictly $\ell^1$-convex and {\rm 2)}~not an $\bR^3$-cross. 
\endproclaim

In other words, for a closed set $M$ its strict $\ell^1$-convexity is necessary
for $M$ to be a sun in~$\ell^\infty(n)$, and is sufficient provided that $M$ is not
a cross (for~$\ell^\infty(3))$ or not a cocross (for $\ell^\infty(n))$.

In the plane $\bR^2$ the situation is trivial: any cocross in~$\bR^2$ 
is an $\bR^2$-cross, which evidently is not a strict  
$\ell^1$-convex set. In follows that in~$\bR^2$ Theorem~A coincides with
Theorem~1.

The next proposition characterises strict $\ell^1$-convex 
cocrosses among all cocrosses in~$\bR^3$. We will prove it later.  

\proclaim{Proposition~1} A cocross in~$\bR^3$ is strictly $\ell^1$-convex
if and only if it is a connected $\bR^3$-cross.
\endproclaim

(Topological dimension of an $\bR^n$-cross $K$ is always~1, therefore its
connectedness is equivalent to arcwise connectedness and, because of 
the inclusion $K\subset\cross(x)$, is equivalent to $\ell^1$-convexity.)

\subhead{3. Preliminary results}
\endsubhead

As usual, by $e_1$, $e_2$, $e_3$ we denote the standard basis of 
unit vectors in~$\bR^3$.

\proclaim{Lemma 1} Let $\emptyset\ne M\subset\bR^3$ has the following properties: 
\roster
\item"{\rm 1)}"~$0$, $e_1$, $e_2$, $e_3\in M$,
\item"{\rm 2)}"~$M$ is not a cross; 
\item"{\rm 3)}"~any two points $x,y\in M$ such that $x\ff y$ can be connected
by a curve $k(t)\subset M$, $k(0)=x$, $k(1)=y$, having strictlty monotonic
coordinate functions. 
\endroster
Then there is a point  $w\in M$ such that  $w\ff0 $.
\endproclaim

\demo{Proof} Let $x\in M$ and let $x\notin \cross(0)$. Such a point exists
by the second assumption. If $x\ff 0$ then the desired statement
is proved. Suppose that $x\in c^3(0)$. Without loss of generality
assume that $x_1=0$. Since $M\not\subset \cross(0)$ then $x\ff e_1$. From
our assumptions it follows that $x$ and~$e_1$ are connected by the curve 
$k(t)\subset M$, $k(0)=x$, $k(1)=e_1$, with strictly monotonic in~$t$ 
coordinate functions. It is clear that $c^3(0)\subset 
c^3(e_1)\cup c^3(x)$. Now the strict monotonicity of functions
$k_i(t)$, $i=1,2,3$, implies that  $k(t)\notin 
c^3(e_1)\cup c^3(x)$ for $0<t<1$, and therefore  
$k(t)\notin c^3(0)$; i.e., $k(t)\ff 0$ for every $t\in (0,1)$. 
This proves the lemma, since $k(t)\subset M$.~$\qed$
\enddemo

For $x\in X$, $P_Mx$ denotes the set of all nearest points from~$M$
to~$X$; i.e., 
$$
  P_Mx=\{y\in M\,|\, \|x-y\|=\rho(x,M):=\inf_{z\in M}\|x-z\|\}.
$$

The next evident lemma is well known (see e.g.\ [9]).

\proclaim{Lemma 2} Let $X$ be a normed linear space 
and let $M\subset X$, $x\in X$, $x\notin M$. Suppose that  
$\hat y\in M$ is a solar point from~$M$ to~$x$ and $y\in P_Mx$. 
Then $[y,\hat y]\subset S(x,\|x-y\|)$.
\endproclaim

The next result will be used in proof of Proposition~1.

\proclaim{Lemma 3} Suppose that a set $M\subset\bR^3$ is  
$\ell^1$-convex and is not a cross. 
Then there exist points $x,y\in M$ such that $x\ff y$.
\endproclaim

\demo{Proof} Suppose that the statement is not true. Without loss 
of generality assume that $0\in M$ and $\eqc(M)=\emptyset$. 
Since $M\not\subset \cross(0)$, 
there exists a point $x\in M$ such that $x\in c^3(0)$, 
$x\notin \cross(0)$. (Note that the case $M\not\subset 
c^3(0)$ immediately implies the desired statement.) 

Without loss of generality we assume $x=(x_1,x_2,0)$, $x_1, x_2\ne 0$. 
Since $\eqc(M)=\emptyset$, there exists a point $y\in M$ such that  
$y_3\ne 0$. By the assumption,  $\eqc(x,y)\ne \emptyset$, therefore  
either $y_1=x_1$, $y_2=0$, or $y_2=x_2$, $y_1=0$. Without loss of 
generality we assume that the second case is fulfilled; i.e., $y=(0,x_2,y_3)$. 
Denote $\xi=(0,x_2,0)$. From $\ell^1$-convexity of~$M$ it follows that
$[y,\xi]\cup [\xi,x]\subset M$. Again applying $\ell^1$-convexity, we have
$[0,\xi]\subset M$.

Now we have $\xi$, $0$, $x$, $y\in M$. After the appropriate translation
of the origin into the point~$\xi$ and the appropriate combination of orthogonal
pivots and extensions, the points  $x,y,0$ will become $e_1, e_2, e_3$ 
in the new coordinate system, the point $\xi$ will be its origin, and the
set $M$ will be represented as~$\tilde M$, having the same properties as
$M$ had. By Lemma~1, there is a point $\tilde w\in \tilde M$, $\tilde w\ff 0$. 
The point  $\tilde w$ in the original system will become $w$, 
for which will be $w\ff \xi$. Now the desired pair $(w,\xi)$, $w\ff \xi$ of
points is found. Lemma~3 is proved.~$\qed$
\enddemo

\subhead{4. Proof of Theorem~1}
\endsubhead

In~[1] it is proved that a strict
sun in~$\ell^\infty(n)$ is strictly $\ell^1$-convex and is not a cocross.
Since a cross is certainly a cocross, we proved the {\smc ``only if''}
part in the theorem. 

{\smc ``If''.} It is enough [1] to consider the following setting: 
$\eqc(M)=\emptyset$.

By Berens--Hetzelt theorem the set $M$ is a sun. We will prove that $M$
is a strict sun by applying the well-known lemma of E.\,B.~Oshman and A.~Br\o{}ndsted
(see e.g.\ [4, Chapter~3]); in accordance with this result it will be enough
to prove that  
$$
  \forall x\notin M\
  \text{and}\ \forall y\in P_Mx\  \Ko(y,x)\cap M=\emptyset,
  \tag\ 1
$$ 
where
$$
  \Ko(y,x)=\{z \mid [z,y]\cap \Bo(x,\|x-y\|)\ne\emptyset\}
$$
is the open supporting cone for the ball $B(x,\|x-y\|)$ at the point $y$ 
(equivalent definitions of supporting cones one can find in [4, Chapter~3]). 

Without loss of generality we assume $x=0$, $\rho(0,M)=\|y\|=1$, and $y_i\ge 0$, 
$i=1,2,3$. Let us denote $J=\{j\,|\, y_j=1\}$. Since $\|y\|=1$ and $y_i\ge 0$, 
we have $J\ne\emptyset$.

Suppose that (1) is false; i.e., $\Ko(y,0)\cap M\ne\emptyset$.

Let $H=c^3(y)\cap\Ko(y,0)$.

1. Suppose at first that $w\in \Ko(y,0)\cap M$ and $w\notin H$. It 
follows that $w\notin c^3(y)$. Therefore $w\ff y$, and, by the
strict $\ell^1$-convexity of $M$, it implies that  $w$ and~$y$ can be 
connected by a strictly monotonic $\ell^1$-geodesic segment $k(t)\subset 
M$, $k(0)=y$, $k(1)=w$.

Let us prove that $k(t)\subset\Bo(0,1)$ for all sufficiently small~$t$.
Under our assumptions we have  $\Ko(y,0)=\{z\,|\, z_j<1$ for $j\in
J\}$. Since $w\in\Ko(y,0)$ and since coordinate functions of the curve
$k(t)$ are strictly monotonic, we have $w_j<k_j(t)<1\quad$ for $j\in J$ and
$0<t<1$.

Since $k_j(t)\to y_j=1$ for $t\to 0$ and $j\in J$, it follows that for
some $\eps>0$
$$
  0<k_j(t)<1 \quad \text{при\ всех}\ j\in J\ \ \text{и}\ \ t\in 
  (0,\eps).  \tag 2
$$
Further, $k_m(t)\to y_m$ as $t\to 0$, $1\le m\le 3$, and $0\le 
y_i<1$ for $i\notin J$. Therefore, for some $\eps_1>0$ we have 
$$
  |k_i(t)-y_i|<1-y_i \quad \text{при}\ \ i\notin J,\ \ 
  0<t<\eps_1.
$$
From this and from (2) it follows that 
$$
   -1 < k_m(t) < 1\quad\text{при}\ \ t\in 
  (0,\min\{\eps,\eps_1\}),\ \
        m=1,2,3.
$$
Therefore $k(t)\subset \Bo(0,1)$ for $t\in (0,\min\{\eps,\eps_1\})$.
But this is impossible, since $y\in P_M 0$. 

2. Suppose now that $M\cap\Ko(y,0)\subset H$. Fix $u\in H\cap M$. 
Since $M$ is stricly $\ell^1$-convex, by Lemma~3 there are points
$\xi,\eta\in M$, $\xi\ff\eta$ which are connected
by a curve $\varkappa(t)\subset M$, $\varkappa(0)=\xi$, $\varkappa(1)=\eta$ having 
stricly monotonic in~$t$ cordinate functions. This implies that  
there is a point $v=\varkappa(t_0)\in M$, $0<t_0<1$, such that $v\ff u$.
From strict $\ell^1$-convexity of~$M$ it follows that $u$ and~$v$ can be 
connected by a strictly monotonic $\ell^1$-geodeisc segment~$k(t)\subset M$, 
$0\le t\le 1$, $k(0)=u$, $k(1)=v$. Since $u\in\Ko(y,0)$, we have that  
$u_j<1$ for $j\in J$ and therefore $k_j(t)<1$ for $j\in J$, 
$0<t<\eps_2$ for some $\eps_2$; i.e, $k(t)\in\Ko(y,0)$ 
for  $0<t<\eps_2$. From strict monotonicity of~$k(t)$ it follows  
that there exists $\eps_3>0$ such that $k(t)\ff y$ for all 
$t\in (0,\eps_3)$. Now we have that for $0<t<\min\{\eps_2,\eps_3\}$ the
curve  $k(t)\subset M$ is contained in $\Ko(y,0)$ and is not contained
in~$H$, which contradicts the assumption $M\cap\Ko(y,0)=H$.
Thus the assumption $\Ko(y,0)\cap M\ne\emptyset$ is false and by 
Oshman--Br\o{}ndsted's lemma $M$ is a strict sun. 
The theorem is proved.~$\qed$
\smallskip

\demo{Proof of Proposition~{\rm 1}} The {\smc ``Only if''} part. 
Let $K\subset\bR^3$ be an $\ell^1$-convex cocross which is not 
an $\bR^3$-cross. Without loss of generality we assume that $K\subset c^3(0)$.

Let us at first consider the case when $\eqc(K)\ne\emptyset$.  We can assume
that  $\{3\}\in\eqc(K)$ and that $K\subset c^2(0)$.
By the definition of a cocross $\card\eqc(K)< 2$, whence  
$\eqc(K)=\{3\}$; i.e., $\eqc_{\bR^2}(K)=\emptyset$. It means that $K$ has
points of the form $x=(\xi,0)$ and $y=(0,\eta)$, $\xi,\eta\ne\emptyset$.
It is clear that  $x\ffm K y$. Hovewer, since $K\subset c^2(0)$, it follows
that $x$ and $y$ cannot be connected by a curve $k(t)\subset 
K$ with strictly monotonic in~$t$ coordinate functions $k_i(t)$, $i=1,2$,
which contradicts the strict $\ell^1$-convexity of~$K$.

\smallskip

Now let us consider the case when $\eqc(K)=\emptyset$. Since $K$ is 
not a cross, by Lemma~3 there are $x,y\in K$ such that $x\ff y$. 
By the strict $\ell^1$-convexity of~$K$ such $x$ and~$y$ are connected 
by a curve $k(t)\subset K$ with strictly monotonic in~$t$ coordinate functions 
$k_i(t)$, $i=1,2,3$. It is easy to see that this is impossible, since 
$K\subset c_{\bR^3}(0)$.

This contradiction shows that a stricty $\ell^1$-convex cocross $K$ has
to be an $\bR^3$-cross, which is $\ell^1$-convex, and therefore connected.
\smallskip

The {\smc ``If''} part. We will prove that every connected  $\bR^n$-cocross $K$
for $n\ge 3$ is strictly $\ell^1$-convex. We may assume $K\subset c^n(0)$. 
By definition, $\eqc(K)=\emptyset$. Every point from~$K$ is of the form: 
$(0,\dots,\xi,0,\dots,0)$. Since $n\ge 3$, we see  that every two points
from~$K$ have at least one common coordinate; i.e.,  $x\ff y$ is never
fulfilled for $x,y\in K$.~$\qed$
\enddemo

\proclaim{Remark} \rm In $\bR^4$ there exists a strictly $\ell^1$-convex
cocross which is not a cross. In fact, let us consider the set 
$$
        K=\bigg\{x=(x_1,x_2,x_3,x_4)\in\bR^4\,\bigg|\,
        \matrix
                x_1=x_2=x_3=0,      & \text{if}\ & x_4\ne 0 \\
                \card\eqc(x,0)\ge 1,& \text{if}\ & x_4 =  0
        \endmatrix\bigg\}.
$$
(The set $K$ is union of the main cocross $c^3(0)$ in~$\bR^3$ 
and of the line $\{x\in\bR^4\,|\, x_1=x_2= x_3=0\}$ which is orthogonal 
to $c^3(0)$ in~$\bR^4$.)  
Now $K\subset c^4(0)$ and since 
$\eqc(K)=\emptyset$, $K$ is a cocross. 
Hovewer, all two points from~$K$ have at least one common coordinate.  
Now the strict $\ell^1$-convexity of~$K$ is evident.  
\endproclaim

\end